\documentclass[11pt]{amsart}

\usepackage{pdfsync}
\usepackage{geometry,graphicx}
\usepackage{amsmath,amssymb,amsthm,esint,bbm}
\usepackage[active]{srcltx}

\usepackage{graphicx}
\catcode`@=11 \@addtoreset{equation}{section} \catcode`@=12

\usepackage{setspace}

\calclayout

\newcommand{\eea}{\end{eqnarray}}

\newcommand{\disp}{\displaystyle}

\newcommand{\beginc} {\begin{center}}

\newtheorem{thm}{\textbf{Theorem}}[section]
\newtheorem{example}{\textbf{Example}}
\newtheorem{lem}[thm]{\textbf{Lemma}}

\newtheorem{rem}[thm]{\textbf{Remark}}
\newtheorem{cor}[thm]{\textbf{Corollary}}
\newtheorem{defn}[thm]{\textbf{Definition}}
\theoremstyle{remark}

\theoremstyle{definition}
\newtheoremstyle{Claim}{}{}{\itshape}{}{\itshape\bfseries}{:}{ }{#1}
\theoremstyle{Claim}

\begin{document}
\title{An improved compact embedding theorem for degenerate Sobolev spaces} 
\author{Dario D. Monticelli and Scott Rodney}
\address{D. D. Monticelli \\
Dept. of Mathematics \\
Politecnico di Milano \\
 20133, Milano Italy}
\email{dario.monticelli@gmail.com}
\address{Scott Rodney\\
Dept. of Mathematics, Physics and Geology \\ 
Cape Breton University \\
Sydney, NS Canada} 
\email{scott\_rodney@cbu.ca}
\thanks{D.~D.~Monticelli is a member of the Gruppo Nazionale per l'Analisi Matematica, la Probabilit\`{a} e le loro Applicazioni (GNAMPA) of the Istituto Nazionale di Alta Matematica (INdAM) and is supported by the PRIN-2015KB9WPT Grant ``Variational methods with applications to problems in mathematical physics and geometry". \\
S.~Rodney is supported by the Canadian NSERC Discovery
  Grant ``Degenerate Elliptic Equations: Regularity of weak solutions with applications".}
\date{March, 2018}

\maketitle

\begin{abstract}
This short note investigates the compact embedding of degenerate matrix weighted Sobolev spaces into weighted Lebesgue spaces. The Sobolev spaces explored are defined as the abstract completion of Lipschitz functions in a bounded domain $\Omega$ with respect to the norm:
$$\|f\|_{QH^{1,p}(v,\mu;\Omega)} = \|f\|_{L^p_v(\Omega)} + \|\nabla f\|_{\mathcal{L}^p_Q(\mu;\Omega)}$$
where the weight $v$ is comparable to a power of the pointwise operator norm of the matrix valued function $Q=Q(x)$ in $\Omega$.  Following our main theorem, we give an explicit application where degeneracy is controlled through an ellipticity condition of the form
$$w(x)|\xi|^p \leq \left(\xi\cdot Q(x)\xi\right)^{p/2}\leq \tau(x)|\xi|^p$$
for a pair of $p$-admissible weights.  We also give explicit examples demonstrating the sharpness of our hypotheses. 
\end{abstract}

\section{Introduction}

  In the study of possibly degenerate elliptic partial differential equations of second order, existence of weak solutions, be they variational or otherwise created, the compactness of embeddings of degenerate matrix weighted Sobolev spaces into Lebesgue spaces plays an important role, see for example \cite{MR} and \cite{CW}.  Our work here improves the local results found in \cite[section 3]{CRW}, specifically \cite[ corollary 3.20]{CRW}. The perspective we take in this work was already thought of in \cite{R} and considers the problem in terms of Sobolev spaces constructed with respect to completions of Lipschitz functions with respect to certain weighted/unweighted Lebesgue norms.  The advantage of this is that we do not require any notion of a Myers-Serrin $H=W$ result.
  
  The plan of this paper is as follows.  In the introduction below, we describe our results in the context of the degenerate Sobolev spaces defined in Definition \ref{spaces}. The principal result of this paper is Theorem \ref{main} with others deduced as corollaries or consequences of slight modification to the proof of Theorem \ref{main}.  The geometric conditions guiding our work together with a careful description of the required Sobolev and Poincar\'e inequalities with related constructions are given in section 2.  Section 3 contains the proofs of our main results.  In section 3 we give an example of our main result applied in the context where degeneracy is encoded by a pair of weights admissible in the sense of \cite{ChW}.  Section 4 contains a counter example demonstrating that in our setting, a global compact embedding may not hold under the hypotheses of local Sobolev and Poincar\'e inequalities.    
  
  Recall that a normed linear space $X$ is compactly embedded in a normed linear space $Y$ if there is a mapping $\mathcal{P}:X\rightarrow Y$ so that given any bounded sequence $\{x_i\}_i\subset X$, the image sequence $\{\mathcal{P}(x_i)\}_i$ contains a convergent subsequence in $Y$.  In the context of studying regularity of weak solutions to possibly non-linear PDEs, the normed linear space $X$ is the possibly degenerate matrix weighted Sobolev space $QH^{1,p}(v,\mu;\Omega)$.  To be precise, let $\Omega$ be a domain (bounded open and connected subset) of $\mathbb{R}^n$ and let $\mu$ be a regular measure on $\Omega$ absolutely continuous with respect to Lebesgue measure.  Fix also a $\mu$-measurable function $Q:\Omega \rightarrow S_n$ (each matrix entry is $\mu$-meas.) taking values in the collection $S_n$ of all non-negative definite self adjoint matrices.  Given any $1\leq p<\infty$, we consider $v_0$ associated to $Q$ given by $v_0(x)=\|Q(x)\|_{\text{op}}^{p/2}$
  where $$\|Q(x)\|_\text{op} = \disp\sup_{|\xi|=1,~\xi\in\mathbb{R}^n}\left|Q(x)\xi\right|$$
  is the operator norm of the matrix $Q(x)$.  We will use at all times the a weight $v$ on $\Omega$ (that is, a nonnegative locally integrable (w.r.t. $\mu$) function on $\Omega$) such that there exists $c_2>0$ so that 
\begin{eqnarray} v_0(x) \leq c_2v(x); ~x\in\Omega.\nonumber
\end{eqnarray}
Consider now the collection $Lip_0(\Omega)$ of all those locally Lipschitz functions with compact support in $\Omega$.  We define a norm on this collection by setting 
\begin{eqnarray}\label{sobnorm}
\|f\|_{QH^{1,p}(v,\mu;\Omega)} &=& \| f \|_{L^p_v(\mu;\Omega)} + \|\sqrt{Q}~\nabla f\|_{L^p(\mu;\Omega)}\\
&=&\left(\int_\Omega \left|f\right|^pv~d\mu\right)^{1/p} + \left(\int_\Omega \left|\sqrt{Q}\nabla f\right|^p~d\mu\right)^{1/p}.\nonumber
\end{eqnarray}
Note that the norm is well defined for any $f\in Lip_0(\Omega)$ since $\sqrt{Q}\in L^p(\mu;\Omega)$.  Note that for such $f$,  $\nabla f$ exists $\mu$-a.e. by the Rademacher-Stepanov theorem.

\begin{defn}\label{spaces} For $1\leq p<\infty$, the Sobolev space $QH^{1,p}_0(v,\mu;\Omega)$ is defined as the completion of $Lip_0(\Omega)$ with respect to the norm \eqref{sobnorm}.  The Sobolev space $QH^{1,p}(v,\mu;\Omega)$ is the completion of $Lip_Q(\Omega)$ with respect to the same norm \eqref{sobnorm} where $Lip_Q(\Omega)$ is the collection of all those locally Lipschitz functions $f$ defined in $\Omega$ for which $\|f\|_{QH^{1,p}(v,\mu;\Omega)}<\infty$.  It is clear that $QH^{1,p}_0(v,\mu;\Omega)\subset QH^{1,p}(v,\mu;\Omega)$.
\end{defn}

\begin{rem}  Note that \eqref{sobnorm} and the semi-norm $$\left(\int_\Omega \left|\sqrt{Q}\nabla f\right|^p~d\mu\right)^{1/p}$$
are not in general equivalent on $QH^{1,p}_0(v,\mu;\Omega)$ unless a global Sobolev-type inequality of the form $$\int_\Omega |f|^p~vd\mu \leq C\int_\Omega \left|\sqrt{Q}\nabla f\right|^p~d\mu$$
holds for every $f\in Lip_0(\Omega)$.
\end{rem}


Although $QH^{1,p}(v,\mu;\Omega)$ is a collection of equivalence classes of $Lip_Q(\Omega)$ sequences Cauchy with respect to the norm \eqref{sobnorm} we will take the sound perspective of \cite{CRW}, \cite{MRW}, \cite{MRW1}, and \cite{SW2} by identifying both $QH^{1,p}(v,\mu;\Omega)$ and $QH^{1,p}_0(v,\mu;\Omega)$ with a closed subset of $L^p_v(\mu;\Omega)\times \mathcal{L}^p_Q(\mu;\Omega)$ using a natural isometry; see \cite{CRW} for more details.

\begin{rem}\label{pairs} The space $\mathcal{L}^p_Q(\mu;\Omega)$ is the collection of all $\mu$-measurable vector valued functions $\vec{k}$ for which
\begin{eqnarray}
\| \vec{k} \|_{\mathcal{L}_Q^p(\mu;\Omega)} = \left(\int_\Omega \left| \sqrt{Q(x)}~\vec{k}(x)\right|^p~ d\mu\right)^{1/p} <\infty\nonumber
\end{eqnarray}
In the case when $\mu$ is Lebesgue measure, completeness of $\mathcal{L}^p_Q(\mu;\Omega)$ is found in \cite{CRR} with the case $p=2$ first treated in \cite{SW2}. When $\mu$ is a general regular measure absolutely continuous with respect to Lebesgue measure, completeness can be shown using the same techniques with minor modifications.
\end{rem}

With Remark \ref{pairs}, for the rest of this article we denote elements of $QH^{1,p}(v,\mu;\Omega)$ as pairs $(f,\vec{g})$ for which there is a sequence $\{f_j\}\in Lip_Q(\Omega)$ so that
$$ f_j\rightarrow f \text{ in } L^p_v(\mu;\Omega),\text{ and } \nabla f_j \rightarrow \vec{g} \text{ in } \mathcal{L}^p_Q(\mu;\Omega).$$

Our main results examine the compactness of the projection mapping $\pi:QH^{1,p}(v,\mu;\Omega) \rightarrow L^r_v(\mu;E)$ defined by
$$\pi\left(\vec{\bf u}\right) = \pi\left((u,\nabla u)\right) = u$$
for any $E\Subset\Omega$.  We now state our main theorem leaving the technical definitions to the next section.  
\begin{thm}\label{main} Let $1< p < \infty$ and $\rho$ be a quasimetric defined in $\Omega$ whose open balls satisfy Definition \ref{doubling}. Then, $QH^{1,p}(v,\mu;\Omega)$ is compactly embedded in $L^q_v(\mu;E)$ for any $E\Subset\Omega$ and $q\in [1,p\sigma)$ provided the pair $(\Omega,\rho)$ satisfies each of the following.
\begin{center}
\begin{enumerate}
\item $(\Omega,\rho)$ admits a local Poincar\'e inequality of order $p$
\item $(\Omega,\rho)$ admits a local Sobolev property of order $p$ and gain factor $\sigma>1$.
\end{enumerate}
\end{center}
\end{thm}

Given any open subdomain $E\Subset \Omega$, extending Lipschitz functions with compact support in $E$ by zero allows one to consider $QH^{1,p}_0(v,\mu;E)$ as a subspace of $QH^{1,p}(v,\mu;\Omega)$.  This leads immediately to the following corollary of Theorem \ref{main}.

\begin{cor}\label{maincor} Let $1< p < \infty$ and $E\Subset \Omega$. Then, $QH^{1,p}_0(v,\mu;E)$ is compactly embedded in $L^q_v(\mu;E)$ for every $q\in [1,p\sigma)$ provided the pair $(\Omega,\rho)$ satisfies each of the following.
\begin{center}

\begin{enumerate}
\item $(\Omega,\rho)$ admits a local Poincar\'e inequality of order $p$.
\item $(\Omega,\rho)$ admits a local Sobolev property of order $p$ and gain factor $\sigma>1$.
\end{enumerate}
\end{center}
\end{cor}

We also mention a result that is a consequence of the proof of Theorem \ref{main}; this will be clearly pointed out in section 3.

\begin{thm}\label{[1,p)}
Assume the hypotheses of Theorem \ref{main} omitting item (2).  Then, $QH^{1,p}(v,\mu;\Omega)$ is compactly embedded in $L^q_v(\mu;\Omega)$ for $1\leq q<p$.  Further, given any $E\Subset\Omega$, $QH^{1,p}(v,\mu;\Omega)$ is compactly embedded in $L^q_v(\mu;E)$ for $1\leq q \leq p$. 
\end{thm}

As a last result we give conditions similar to those of Theorem \ref{main} under which one obtains compact embedding on all of $\Omega$, not only on $E\Subset \Omega$.

\begin{thm}  Let the hypotheses of Theorem \ref{main} hold with the exception that we replace item (2) with
\begin{itemize}
\item[(2*)] $(\Omega,\rho)$ admits a global Sobolev property of order $p$ and gain $\sigma>1$; see Definition \ref{globsob}. 
\end{itemize}
Then both $QH^{1,p}_0(v,\mu;\Omega)$ and $QH^{1,p}(v,\mu;\Omega)$ are compactly embedded in $L^q_v(\mu;\Omega)$ for any $q\in [1,p\sigma)$.
\end{thm}
While we will not prove this theorem explicitly, the result is gleaned from Remark \ref{Omega} and an interpolation inequality 
$$\|g\|_{L^q_v(\mu;\Omega)}\leq \|g\|_{L^1_v(\mu;\Omega)}^\lambda\|g\|_{L^{p\sigma}_v(\mu;\Omega)}^{1-\lambda}$$
valid for any $g\in L^{p\sigma}_v(\mu;\Omega)$.

In section 4 we present an application of our results to degenerate Sobolev spaces where degeneracy is controlled by $p$-admissible weights in $\Omega$.  See section 4 for complete details.

\begin{thm}\label{2weight}
Fix a bounded domain $\Omega$ of $\mathbb{R}^n$.  Let $w\leq \tau$ be a pair of $p$-admissible weights, for some $1<p<+\infty$ and let $Q(x)$ be a non-negative definite matrix function that satisfies the ellipticity condition
$$w|\xi|^p \leq \left|\sqrt{Q(x)}\xi\right|^p \leq \tau|\xi|^p.$$
Then, there is a $q>p$ so that $QH^{1,p}_0(\tau,dx;E)$ is compactly embedded in $L^r_\tau(E)$ for all $1\leq r<q$ and $E\Subset\Omega$.  Further, $QH^{1,p}(\tau,dx;\Omega)$ is compactly embedded in $L^r_\tau(E)$ for $1\leq r<q$ and any $E\Subset \Omega$.
\end{thm}

\section{Preliminaries}

We begin this section by recalling the quasimetric structure upon which our result is built.  We assume there is a quasimetric $\rho$ on $\Omega$.  That is, there is a $\kappa\geq 1$ so that for each $x,y,z\in\Omega$
\begin{eqnarray}
(i)&&\rho(x,y) \geq 0 \text{ with equality only if $x=y$}\nonumber\\
(ii)&&\rho(x,y) = \rho(y,x)\nonumber\\
(iii)&&\rho(x,y) \leq \kappa\left(\rho(x,z) + \rho(z,y)\right).
\end{eqnarray}
Given $x\in \Omega$ and $r>0$, we denote the $\rho$-ball centered at $x$ with radius $r$ by $$B(x,r) = \{y\in \Omega~:~ \rho(x,y)<r\}.$$
We require that $\rho$-balls are open which is equivalent to the condition 
$$\displaystyle\lim_{y\rightarrow x} \rho(x,y) = 0$$
for every $x\in \Omega$.  Moreover, we require that for every $x\in\Omega$ there is a $\delta=\delta(x)>0$ so that 
$$\overline{B(x,r)} \subset \Omega$$
for any $0<r<\delta$.  As in \cite{CRW}, we will not need to assume that the family $\{B(x,r)\}_{x\in\Omega,r>0}$ admits a doubling measure but we do require a local geometric doubling condition.  
\begin{defn} \label{doubling}
A quasimetric space $(\Omega,\rho)$ is locally geometrically doubling if given any compact subset $K$ of $\Omega$, there is a $\delta>0$ so that $0<s\leq r<\delta$ and $x\in K$ imply that $B(x,r)$ may contain at most $C(r/s)$ centers of disjoint $\rho$-balls of radius $s$; here $C:(0,\infty)\rightarrow (0,\infty)$ is independent of $K$.
\end{defn}
\begin{rem} This condition is weaker than the existence of a locally doubling measure for $\rho$-balls.  We refer the reader to \cite{HyM} for further details and discussions.
\end{rem}

The local geometric doubling condition is used to establish the following lemma giving coverings of compact sets by $\rho$-balls with finite overlaps.  We omit the proof here and point the reader to the proof of \cite[lemma 3.12]{CRW}.

\begin{lem}\label{finiteoverlap}
Let $K$ be a compact subset of $\Omega$ and $c_0\geq 1$.  Then, there are positive constants $\delta_0=\delta_0(K,\kappa,c_0)$ and $P=P(\kappa,c_0)$ so that for any $0<r<\delta_0$ there is a finite collection of $\rho$-balls $\{B(x_j,r)\}_{j=1}^N$, each centred in $K$, that satisfies
    \begin{enumerate}
    \item[(i)] $K\subset \displaystyle\bigcup_{j=1}^N B(x_j,r)\subset \displaystyle\bigcup_{j=1}^N \overline{B(x_j,c_0r)}\subset \Omega$
    \item[(ii)] $\displaystyle\sum_{j=1}^N \chi_{B(x_j,c_0r)}(x) \leq P$ for any $x\in \displaystyle\bigcup_{j=1}^N B(x_j,r)$
    \end{enumerate}
\end{lem}

Associated to our collections of $\rho$-balls are the Sobolev and Poincar\'e inequalities that form the main hypotheses of Theorem \ref{main}.  

\begin{defn}\label{pc} We say that $(\Omega,\rho)$ supports a local Poincar\'e property of order $p$ if there is a $c_0\geq 1$ such that given any compact subset $K$ of $\Omega$ and $\epsilon>0$, there is a $\delta_1>0$ so that $0<r<\delta_1$ and $x\in K$ give
\begin{eqnarray}\label{PE}\|f-f_{B(x,r)}\|_{L^p_v(\mu;B(x,r))} < \epsilon\| (f,\nabla f)\|_{QH^{1,p}(v,\mu;B(x,c_0r))}
\end{eqnarray}
for any $f\in Lip_{loc}(\Omega)$, where $f_B=\frac{1}{v(B)}\int_Bfv\,d\mu$.
\end{defn}

\begin{rem}\label{prem}
Inequality \eqref{PE} may feel unfamiliar.  The reader may be more familiar with the standard $Q$-weighted Poincar\'e inequality:
$$\left(\frac{1}{v(B)}\int_B |f-f_B|^{p}v~d\mu\right)^{1/p} \leq Cr\left(\frac{1}{\mu(c_0B)}\int_{c_0B} \left|\sqrt{Q}\nabla f\right|^p~d\mu\right)^{1/p}$$
holding for any $f\in Lip(B)$.  It is not difficult to see that this inequality is enough to ensure Definition \ref{pc} provided 
\begin{eqnarray}\label{compat}
\disp\lim_{r\rightarrow 0}\sup_{x\in\Omega} \left[r^p\frac{v(B(x,r))}{\mu(B(x,C_0r))}\right] &=&0.
\end{eqnarray}
\end{rem}

\begin{defn}\label{sobolev} We say that $(\Omega,\rho)$ supports a local Sobolev property of order $p$ and gain $\sigma$ if there is a $\sigma\geq 1$ so that given any compact $K\subset \Omega$ one can choose $\delta_2>0$ with the property that if $B=B(x,r)$ is centred in $K$ and of radius $0<r<\delta_2$ then
\begin{eqnarray} \label{sob}
\left(\int_{B} \left| f\right| ^{p\sigma}v~d\mu\right)^{1/p\sigma} \leq C(B)\|(f,\nabla f)\|_{QH^{1,p}(v,\mu;\Omega)}
\end{eqnarray}
for every $f\in Lip_0(B)$.
\end{defn}

\begin{defn}\label{globsob} We say that $(\Omega,\rho)$ admits a global Sobolev inequality if there is a constant $C>0$ so that 
\begin{eqnarray}\label{globsobineq}
\left(\int_\Omega \left|f\right|^{p\sigma}v~d\mu\right)^{1/p\sigma} \leq C\left[\|f\|_{L^p_v(\mu;\Omega)}+\|\sqrt{Q}\nabla f\|_{L^p(\mu;\Omega)}\right]
\end{eqnarray}
for every $f\in Lip_0(\Omega)$.
\end{defn}

\section{Proofs}
\subsection{Proof of Theorem \ref{main}}

We consider first the case $q=p$.  Fix an open set $E\Subset\Omega$.  Let $\{\vec{\bf u}_n\}_n=\{(u_n,\nabla u_n)\}_n$ be a bounded sequence in $QH^{1,p}(v,\mu;\Omega)$ with upper bound $M$ and let $\epsilon>0$.   Given $0<r<\delta=\min\{\delta_0,\delta_1\}$, Lemma \ref{finiteoverlap} provides a finite collection of $\rho$-balls $\{B(x_j,r)\}_{j=1}^N$ satisfying (i) and (ii) of Lemma \ref{finiteoverlap}.  Further, for each $1\leq j\leq N$, we have 
$$\|f-f_{B(x,r)}\|_{L^p_v(\mu;B(x,r))} < \epsilon\| (f,\nabla f)\|_{QH^{1,p}(v,\mu;B(x,c_0r))}$$
for any $f\in Lip_\text{loc}(\Omega)$ by Definition \ref{pc}.\\

In order to show $\{u_n\}$ is Cauchy in $L^p_v(\mu;E)$ we estimate
\begin{eqnarray}\label{A1}
\displaystyle\sum_{j=1}^N\int_{B_j} \left|u_m-u_n\right|^pv~d\mu &\leq& C_p\displaystyle\sum_{j=1}^N\Big[\int_{B_j}\left|u_m-u_n-(u_n-u_m)_{B_j}\right|^pv~d\mu\\
&&\quad\quad\quad\quad\quad + \left|(u_m-u_n)_{B_j}\right|^pv(B_j)\Big]\nonumber\\
&=& C_p\left(I + II\right)\nonumber
\end{eqnarray}
where $B_j = B(x_j,r)$ and, for an integrable function $g$, $g_B = \fint_B gv~d\mu$ is the $v$-average of $g$.  We estimate $I$ and $II$ separately using different techniques.  Beginning with $I$, we assume that $r<\delta=\delta(E)$ and apply the Poincar\'e inequality \eqref{PE} to find

\begin{eqnarray}\label{partI}
I &\leq &\epsilon^p \displaystyle\sum_{j=1}^N\|(u_n-u_m,\nabla(u_m-u_n))\|_{QH^{1,p}(v,\mu;{B}(x_j,c_0r))}^p\nonumber\\
&\leq& \epsilon^pP \|(u_n-u_m,\nabla(u_m-u_n))\|_{QH^{1,p}(v,\mu;\Omega)}^p\\
&\leq& 2^pPM^p\epsilon^p\nonumber
\end{eqnarray}
where $P$ is the overlap constant for our collection as in Lemma \ref{finiteoverlap}.\\

To estimate item $II$ we use weak convergence.  Indeed, since $\{u_n\}$ is a bounded sequence in $L^p_v(\mu;\Omega)$, it admits a weakly convergent subsequence that we denote by $\{u_n\}$ to preserve the index.  As $v(B_j)$ is finite, the characteristic function 
$\chi_{B_j}(x) \in L^{p'}_v(\mu;{\Omega})$ for every $j$.  Thus, there is $T\in\mathbb{N}$ so that $m,n\geq T$ gives
\begin{eqnarray}\label{partII}
II\leq \displaystyle\sum_{j=1}^Nv^{1-p}(B_j)\left|\int_\Omega \chi_{B_j}(u_m-u_n)v~d\mu\right|^p <\epsilon^p
\end{eqnarray}
Combining \eqref{partI} and \eqref{partII} with \eqref{A1} we find 
\begin{eqnarray}\label{theend}
\|u_m-u_n\|_{L^p_v(\mu;E)}^p <C \epsilon^p(1+2^pM^pP)
\end{eqnarray}
when $m,n\geq T$.  This establishes convergence of our subsequence in $L^p_v(\mu;E)$ and, by H\"older's inequality, also in $L^q_v(\mu;E)$ for $1\leq q \leq p$.  This establishes Theorem \ref{main} for the range $1\leq q   \leq p$.  

\begin{rem}\label{Omega} It is not difficult to now show that  $QH^{1,p}(v,\mu;\Omega)$ is compactly embedded in $L^q_v(\mu;\Omega)$ for the range $1\leq q<p$.  Indeed, fix $\eta>0$ and assume that our set $E$ satisfies $v(\Omega\setminus E)<\eta$.  Then, for our subsequence $\{u_n\}$ constructed above, H\"older's inequality and boundedness in $QH^{1,p}(v,\mu;\Omega)$ show that for any $j,k\in\mathbb{N}$,
\begin{eqnarray}\label{diag}
\|u_j-u_k\|_{L^1_v(\mu;\Omega)} &=& \|u_j-u_k\|_{L^1_v(\mu;E)}+\|u_j-u_k\|_{L^1_v(\mu;\Omega\setminus E)} \nonumber\\
&\leq& \|u_j-u_k\|_{L^p_v(\mu;E)}v(E)^{1/p'} + \|u_j-u_k\|_{L^{p}_v(\mu;\Omega\setminus E)}v(\Omega\setminus E)^{1/p'}\nonumber\\
&\leq& \|u_j-u_k\|_{L^p_v(\mu;E)}v(E)^{1/p'} + 2M\eta^{1/p'}.\nonumber\\
\end{eqnarray}
Since $\{u_j\}$ is Cauchy in $L^p_v(\mu;E)$, choosing $j,k$  sufficiently large shows that $\{u_n\}$ is Cauchy in $L^1_v(\mu;\Omega)$.  That is, we have shown $QH^{1,p}(v,\mu;\Omega)$ is compactly embedded in $L^1_v(\mu;\Omega)$.  Interestingly, we also conclude the same for $L^q_v(\mu;\Omega)$ when $1\leq q<p$ through an appeal to H\"older's inequality. Given $1<q<p$, we may choose $\lambda\in(0,1)$  so that 
\begin{eqnarray}\label{interpolato}
\|u_j-u_k\|_{L^q_v(\mu;\Omega)} &\leq& \|u_j-u_k\|_{L^1_v(\mu;\Omega)}^\lambda\|u_j-u_k\|_{L^{p}_v(\mu;\Omega)}^{1-\lambda}\nonumber\\
&\leq& (2M)^{1-\lambda}\|u_j-u_k\|_{L^1_v(\mu;\Omega)}^\lambda.
\end{eqnarray}
This argument completes the proof of Theorem \ref{[1,p)}.
\end{rem}

We now turn our attention to the range $p<q<p\sigma$.  With $0<r<\delta$ as above, cover $E$ with Euclidean balls $D(x,s)$ where $s=s(x)$ is chosen so that $\overline{D(x,s)}\subset B$. By compactness, we may select $\{D(x_j,s_j)\}_{j=1}^{N_1}$ that covers $E$.  Let $\{\varphi_j\}$ be a partition of unity subordinate to this cover and let $f\in Lip_{loc}(\Omega)$.  With $D(x_j,s_j)\subset B_j=B(x_j,r)$ we see from the Sobolev inequality \eqref{sob} that

\begin{eqnarray}
\int_E \left|f\right|^{p\sigma}v~d\mu &\leq& C_p \displaystyle \sum_j \int_{B_j} |f\varphi_j|^{p\sigma}v~d\mu\nonumber\\
&\leq& \left(C_p\displaystyle \sum_jC(B_j)\left[ \int_{B_j} \left|f\right|^{p}v~d\mu + \int_{B_j}\left|\sqrt{Q}\nabla(f\varphi_j)\right|^p d\mu  \right]\right)^\sigma
\end{eqnarray}
since $0\leq \varphi_j(x)\leq 1$ for each $j$.  The second term splits with integrand bounded above by 
$$C_*\left[\left|\sqrt{Q}\nabla f\right|^p + \left|f\right|^pv\right]$$
since $v \geq c_2^{-1} \|Q(x)\|^{p/2}_{op}$ and where $C_*$ is a constant independent of $f$.  Since the sum is finite, we find a constant $\tilde{C}=\tilde{C}(E)$ so that
\begin{eqnarray}\label{noncompactsobolev}
\|f\|_{L^{p\sigma}_v(\mu;E)} \leq \tilde{C}\|(f,\nabla f)\|_{QH^{1,p}(v,\mu;\Omega)}
\end{eqnarray}
for any $f\in Lip_{loc}(\Omega)$; by density also for any pair $(g,\nabla g)\in QH^{1,p}(v,\mu;\Omega)$. Thus, boundedness in $L^{p\sigma}_v(\mu;E)$ of our sequence $\{u_n\}$ is established.  Given $p<q<p\sigma$, we may choose $\lambda\in (0,1)$ so that
$$\|u_j-u_k\|_{L^q_v(\mu;E)} \leq \|u_j-u_k\|_{L^p_v(\mu;E)}^\lambda \|u_j-u_k\|_{L^{p\sigma}_v(\mu;E)}^{1-\lambda}\leq C(2M)^{1-\lambda}\|u_j-u_k\|_{L^p_v(\mu;E)}^\lambda$$
and we conclude that $\{u_n\}$ is Cauchy in $L^q_v(\mu;E)$.  This completes the proof of Theorem \ref{main}.
\begin{flushright}
$\Box$
\end{flushright}

\subsection{Proof of Corollary \ref{maincor}}

Let $E\Subset\Omega$ and $\{\vec{\bf u}_n\}=\{u_n,\nabla u_n\}$ be a bounded sequence in $QH^{1,p}_0(v,\mu;E)$.  Since each element $\vec{\bf u}_n$ may be viewed as an equivalence class of Cauchy sequences of $Lip_0(E)$ functions, we may choose a representative sequence $\{g_m^n\}_m\subset Lip_0(E)$ converging to $u_n$ in $QH^{1,p}(v,\mu;E)$ norm.  For each $m$, set
\begin{eqnarray}
G_m^n=\left\{\begin{array}{ccc}
 g_m^n(x)&\textrm{ if }&x\in E\\
0 &\textrm{ if }&x\in \Omega\setminus E
\end{array}\right.
\end{eqnarray}

The resulting sequence of extended functions $\{G_m^n\}$ is Cauchy in $QH^{1,p}(v,\mu;\Omega)$ and converges to $\vec{\bf w}_n=(w_n,\nabla w_n)\in QH^{1,p}(v,\mu;\Omega)$ with 
$$\|u_n-w_n\|_{L^p_v(\mu;E)} = \|\nabla u_n-\nabla w_n\|_{\mathcal{L}^p_Q(\mu;E)} =\|\vec{\bf u}_n-\vec{ \bf w}_n\|_{QH^{1,p}(v,\mu;E)} =0.$$
From this we can also see that $u_n=w_n$ in $L^{p\sigma}_v(\mu;E)$.  Since our new sequence $\{\vec{\bf w}_n\}$ is bounded in $QH^{1,p}(v,\mu;\Omega)$, Theorem \ref{main} provides a subsequence of $\{w_n\}$ (that we refer to as $\{w_n\}$ to preserve the index) that is Cauchy in $L^q_v(\mu;E)$ for each $q\in[1,p\sigma)$.  Since 
$$\|u_j-u_k\|_{L^q_v(\mu;E)} = \|w_j-w_k\|_{L^q_v(\mu;E)}$$
for every $j,k$ and $q\in[1,p\sigma]$, we find $\{u_n\}$ is Cauchy in $L^q_v(\mu;E)$ for $1\leq q<p\sigma$.  We now conclude that $QH^{1,p}_0(v,\mu;E)$ is compactly embedded in $L^q_v(\mu;E)$ for each $1\leq q<p\sigma$ completing the proof of Corollary \ref{maincor}.
\begin{flushright}
$\Box$
\end{flushright}



\section{Application to Two Weight Degenerate Problems}

As an application to Theorem \ref{main}, we present compact embeddings for Sobolev spaces with degeneracy controlled by admissible weights.  Given $p> 1$, two weights $w\leq \tau$ on $\Omega$ are called $p$-admissible in $\Omega$ if each of the following conditions are met.
\begin{enumerate}
\item $\tau$ is doubling for the collection of Euclidean balls with center in $\Omega$.  That is, there is a constant $C$ so that given $x\in \Omega$ and $r>0$,
$$\tau(D(x,2r))=\int_{D(x,2r)}\tau~dz \leq C\int_{D(x,r)}\tau~dz=C\tau(D(x,r))$$
\item $w\in A_p(\Omega)$.  For $1<p<\infty$, the Muckenhoupt class of weights $A_p(\Omega)$ is the collection of all those non-negative functions $\varphi\in L^1_\textrm{loc}(\Omega)$ for which
$$\displaystyle\sup_D\left(\frac{1}{|D|}\int_D \varphi~dz\right)\left(\frac{1}{|D|}\int_B \varphi^\frac{1}{1-p}~dz\right)^{p-1}<\infty$$
where the supremum is taken over all Euclidean balls $D=D(x_0,r) = \{x\in\Omega\;:\; |x-x_0|<r\}$ centred in $\Omega$.
\item $w,\tau$ satisfy the Chanillo-Wheeden balance condition; see \cite{ChW} and \cite{CMN}.  That is there are $C>0$ and $q>p$ so that for $0<s\leq r$ and $x\in \Omega$, 
\begin{eqnarray}\label{balancecondition}\frac{s}{r}\left(\frac{\tau(D(x,s))}{\tau(D(x,r))}\right)^{1/q} \leq C \left(\frac{w(D(x,s)}{w(D(x,r)}\right)^{1/p}\end{eqnarray}
\end{enumerate}
where here for a weight $\nu$, $\nu(D) = \int_D \nu(x)~dx$.  For the reader unfamiliar with such objects, power weights $\tau(x)=|x|^t$ form an excellent example category; \cite{CMN} is also a good reference for this deep subject. \cite{CMN} draws from \cite{CW} and other classic works in the area to demonstrate that for admissible weights $w\leq \tau$ there are constants $C>0, q>p$ so that for any Euclidean ball $D=D(x,r)\Subset \Omega$ one has

\begin{enumerate} 
\item the local Poincar\'e inequality
$$\left(\frac{1}{\tau(D)}\int_D |f-f_{D;\tau}|^q\tau~dx\right)^{1/q}\leq 
Cr \left(\frac{1}{w(D)}\int_D \left|\nabla f\right|^pw~dx\right)^{1/p}$$
for any $f\in Lip_\text{loc}(\Omega)$, and
\item the local Sobolev inequality 
$$\left(\frac{1}{\tau(D)}\int_D |g|^q\tau~dx\right)^{1/q}\leq 
C r\left(\frac{1}{w(D)}\int_D \left|\nabla g\right|^pw~dx\right)^{1/p}$$
for each $g\in Lip_0(D)$.
\end{enumerate}
These inequalities are used to study second order degenerate elliptic problems $\left( Xu = \varphi\right)$ with principal part of $X$ given by a matrix weighted $p$-Laplacian: 
$$Lu=\text{Div}\left(\left|\sqrt{Q}\nabla u\right|^{p-2}{Q}\nabla u\right).$$
The symmetric non-negative definite $n\times n$ matrix $Q$ is assumed to satisfy the degenerate elliptic condition
\begin{eqnarray} \label{elliptic}w(x)\left|\xi\right|^p \leq  \left|\sqrt{Q(x)}\xi\right|^{p}\leq v_0(x)\left|\xi\right|^p \leq \tau(x)\left|\xi\right|^p,~ a.e.~ x\in\Omega, \xi\in\mathbb{R}^n
\end{eqnarray}
where $v_0 = \|{Q}\|_{op}^{p/2}$ is the $p^\text{th}$ power of the operator norm of $\sqrt{Q(x)}$.  Weak solution spaces for Dirichlet and Neumann problems associated to such equations are the matrix weighted Sobolev spaces $QH^{1,p}_0(\tau,dx;\Omega)$ and $QH^{1,p}(\tau,dx;\Omega)$, as defined in \S1 with $v\equiv\tau$.  Because of the ellipticity condition \eqref{elliptic}, we find the Sobolev and Poincar\'e inequalities \eqref{sob} and \eqref{PE} of Definitions \ref{sobolev} and \ref{pc}. Indeed, to see that the Poincar\'e holds, let $f\in Lip_\text{loc}(\Omega)$, $E$ a compact subset of $\Omega$ and fix a Euclidean ball $D=D(x,s)$ with $s<r_0=dist(E,\partial\Omega)$.  Using the two weight Poincar\'e estimate we see
\begin{eqnarray}
\left(\int_D \left|f-f_D\right|^p\tau~dx\right)^{1/p} &\leq&
\tau(D)^{\frac{1}{p}}\left(\frac{1}{\tau(D)}\int_D |f-f_D|^q\tau~dx\right)^{1/q}\nonumber\\
&\leq&
Cs\frac{\tau(D)^{\frac{1}{p}}}{w(D)^{\frac{1}{p}}}\left(\int_D|\sqrt{Q}\nabla f|^p~dx\right)^{1/p}.\nonumber
\end{eqnarray}
The balance condition \eqref{balancecondition} with $r=r_0=\text{dist}(E,\partial\Omega)$ gives a positive constant $C_1$ so that 
$$s \left[\frac{\tau(D)^{1/q}}{w(D)^{1/p}}\right] \leq C_1$$
for every $x\in E$.  As a result, we see that 
$$\lim_{s\rightarrow 0}\sup_{x\in E} \left[Cs\frac{\tau(D)^{\frac{1}{p}}}{w(D)^{\frac{1}{p}}}\right] \leq C\lim_{s\rightarrow 0}\sup_{x\in E}\tau(D(x,s))^\frac{q-p}{qp} = 0 \text{ as }q>p$$
and we conclude that Definition \ref{pc} holds for Euclidean balls $D(x,r)$ with $c_0=1$.  The argument giving Definition \ref{sobolev} is similar and left to the reader.  This concludes the proof of Theorem \ref{2weight}.

\section{Failure of Global Compact Embedding Under Local Hypotheses}

In this section we give an example of a matrix function $Q$ satisfying the conditions of Theorem \ref{main} on a specific domain of $\mathbb{R}^n$ where the embedding
$$\pi:QH^{1,p}_0(v,\mu;\Omega)\subset QH^{1,p}(v,\mu;\Omega) \rightarrow L^q_v(\mu;\Omega)$$
fails to be compact for any $q>p>1$, where $v\equiv 1$ and $d\mu=dx$ in $\Omega$.  In particular, there is no embedding at all. A posteriori, this is because no global Sobolev inequality with gain $\sigma>1$ holds on $\Omega$.

\begin{example}
Let $n\geq 2$ and set $\Omega=(0,1)\times...\times (0,1)=(0,1)^n$; the $n$-dimensional unit cube.  Let $p>1$, $q>p$ and choose $\beta\in [\frac{1}{q},\frac{1}{p})$. Define 
\begin{eqnarray}
Q(x) = \textrm{Diag}\left[x_1^2,1...,1\right].\nonumber
\end{eqnarray}
Since $Q(x)$ is uniformly elliptic away from the boundary $\partial\Omega$, standard local Sobolev and Poincar\'e inequalities of order $p$ hold on Euclidean balls contained in $\Omega$, both with gain $\sigma=\frac{n}{n-p}$ when $n>p$ or any $\sigma>1$ if $n\leq p$.  \\

Consider the function $u:\Omega\rightarrow \mathbb{R}^n$ defined by 
$$u(x)  = \left\{
\begin{array}{ccc}(x_1^{-\beta}-2)\psi(\hat{x})&\textrm{ if }&0<x_1<2^{-\frac{1}{\beta}},\\
0&\textrm{ if }&2^{-\frac{1}{\beta}}\leq x_1<1
\end{array}\right.
$$
with gradient
$$\nabla u(x) = \left\{
\begin{array}{ccc}
(-\beta x_1^{-\beta-1}\psi(\hat{x}),(x_1^{-\beta}-2)\nabla_{\hat{x}}\psi(\hat{x}))& \textrm{ if }&0<x_1<2^{-\frac{1}{\beta}},\\
0&\textrm{ if }& 2^{-\frac{1}{\beta}}\leq x_1<1
\end{array}
\right.$$
where $\hat{x}=(x_2,...,x_n)\in(0,1)^{n-1}$ and $\psi\in C^\infty_0((0,1)^{n-1})$. 

From these definitions it is easy to check that $u\in L^p(\Omega)$, $\left|\sqrt{Q}\nabla u\right|\in L^p(\Omega)$ and that $u\notin L^q(\Omega)$. We will now demonstrate a sequence of Lipschitz function with compact support in $\Omega$ that converge to the pair $(u,\nabla u)$ in the $QH^{1,p}(\Omega)$ norm thus showing that $(u,\nabla u)\in QH^{1,p}_0(\Omega)\subset QH^{1,p}(\Omega)$.  

For $j\in\mathbb{N}$, define the $Lip_0(\Omega)$ function (see Figure \ref{fig:AA} for $u_j(t)$ in dimension $n=1$)

\begin{figure}[h]
\includegraphics[width=9cm,scale=1.1]{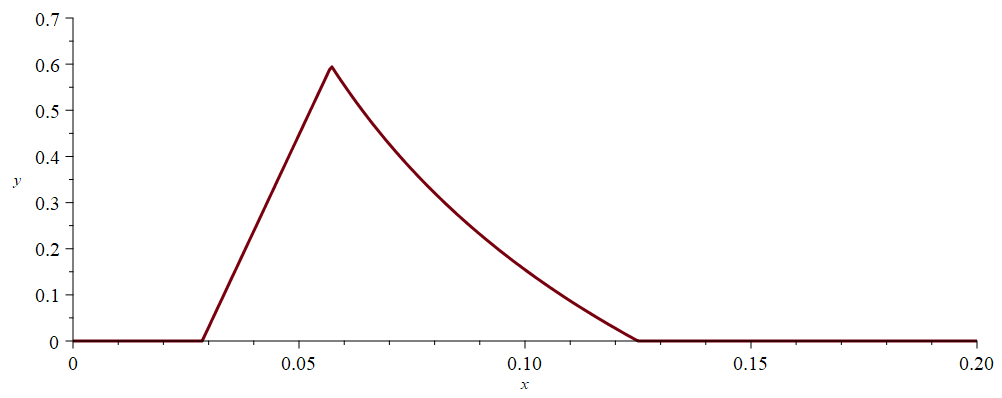}
\caption{$u_j(t)$}
\label{fig:AA}
\end{figure}

\begin{eqnarray}
u_j(x) = \left\{\begin{array}{ccc}
0&\textrm{ if }&0< x_1\leq \frac{1}{j},\\
\left[\left(\left(\frac{j}{2}\right)^\beta-2\right) j(x_1-\frac{1}{j})\right]\psi(\hat{x})&\textrm{ if }& \frac{1}{j}<x_1<\frac{2}{j},\\
(x_1^{-\beta}-2)\psi(\hat{x})&\textrm{ if }&\frac{2}{j}\leq x_1<2^{-\frac{1}{\beta}},\\
0&\textrm{ if }&2^{-\frac{1}{\beta}}\leq x_1<1.
\end{array}
\right.
\end{eqnarray}
defined for $j>2^{1+\frac{1}{\beta}}$.

It is not a difficult exercise to show that $|u-u_j|$ and $\left|\sqrt{Q}\nabla (u- u_j)\right|$ converge to $0$ in $L^p(\Omega)$.  As a result, $\{u_j\}$ is a Cauchy sequence of $Lip_0(\Omega)$ functions in $QH^{1,p}(\Omega)$ norm whose limit is $(u,\nabla u)$. Thus, $QH^{1,p}_0(\Omega)\not\subset L^q(\Omega)$. That is, the embedding fails and is thus, obviously, not compact.
\end{example}

Our next example demonstrates that the lack of a global Sobolev inequality a posteriori causes failure of compact embedding of $QH^{1,p}(\Omega)$ and of $QH^{1,p}_0(\Omega)$  in $L^p(\Omega)$, with $\Omega$ as above.  In order to describe this precisely we first appeal to a one dimensional version.

\begin{figure}[h]
\includegraphics[width=8cm]{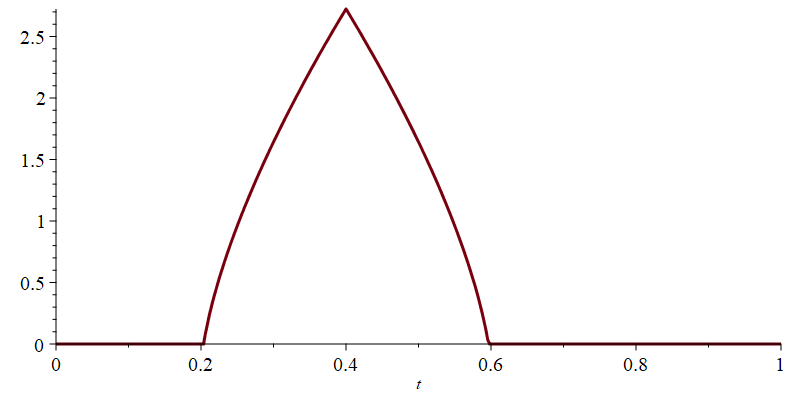}
\caption{$v_j(t)$}
\label{fig:BB}
\end{figure}

\begin{example}
Fix $I=[0,1]$, $p>1$ and set $q(t)=t^{2p}$.  For large $j\in \mathbb{N}$ define the Lipschitz fucntion with compact support in $(0,1)$ (see Figure \ref{fig:BB})
\begin{eqnarray}
v_j(t) =\left\{\begin{array}{ccc}
0&\textrm{ if }&0<t\leq \frac{1}{n}+\frac{1}{n^{p+2}},\\
&&\\
n^\frac{2}{p}\left(t-\frac{1}{n}\right)^\frac{1}{p}-\frac{1}{n}&\textrm{ if }&\frac{1}{n}+\frac{1}{n^{p+2}}<t\leq \frac{2}{n},\\
&&\\
n^\frac{2}{p}\left(\frac{3}{n}-t\right)^\frac{1}{p}-\frac{1}{n}&\textrm{ if }&\frac{2}{n}<t\leq \frac{3}{n}-\frac{1}{n^{p+2}},\\
&&\\
0&\textrm{ if }& \frac{3}{n}-\frac{1}{n^{p+2}}<t\leq 1
\end{array}
\right.
\end{eqnarray}
\end{example}

whose derivative is given by a.e. by

\begin{eqnarray}
v_j'(t)=\left\{\begin{array}{ccc}
0&\textrm{ if }&0<t\leq \frac{1}{n}+\frac{1}{n^{p+2}},\\
&&\\
\frac{1}{p}n^\frac{2}{p}\left(t-\frac{1}{n}\right)^{\frac{1}{p}-1}&\textrm{ if }&\frac{1}{n}+\frac{1}{n^{p+2}}<t\leq \frac{2}{n},\\
&&\\
-\frac{1}{p}n^\frac{2}{p}\left(\frac{3}{n}-t\right)^{\frac{1}{p}-1}&\textrm{ if }&\frac{2}{n}<t\leq \frac{3}{n}-\frac{1}{n^{p+2}},\\
&&\\
0&\textrm{ if }& \frac{3}{n}-\frac{1}{n^{p+2}}<t\leq 1.
\end{array}
\right.
\end{eqnarray}

Clearly, for each $j$ large enough 
$$\int_0^1 v_j(t)^p~dt\leq 1,$$
and using Jensen's inequality and that $v_j(t)$ is concave in the interval $(\frac{1}{n}+\frac{1}{n^{p+2}},\frac{2}{n})$ we see that 
\begin{eqnarray}\int_0^1 v_j^p(t)~dt &=&
2\int_{\frac{1}{n}+\frac{1}{n^{p+2}}}^\frac{2}{n} \left(n^\frac{2}{p}\left(t-\frac{1}{n}\right)^\frac{1}{p}-\frac{1}{n}\right)^p~dt\nonumber\\
&\geq& 2\left(\int_{\frac{1}{n}+\frac{1}{n^{p+2}}}^\frac{2}{n} \left(n^\frac{2}{p}\left(t-\frac{1}{n}\right)^\frac{1}{p}-\frac{1}{n}\right)~dt\right)^p\left(\frac{1}{n}-\frac{1}{n^{p+2}}\right)^{1-p}\nonumber\\
&\geq& 2\left(\frac{1}{2}\left(\frac{1}{n}-\frac{1}{n^{p+2}}\right)\left(n^\frac{1}{p}-\frac{1}{n}\right)\right)^p\left(\frac{1}{n}-\frac{1}{n^{p+2}}\right)^{1-p}\nonumber\\
&\rightarrow& 2^{1-p} ~~~\textrm{ as } j\rightarrow\infty.\nonumber
\end{eqnarray}
Then, the sequence $\{v_j\}$ is bounded in $L^p(I)$, it converges to zero pointwise everywhere in $I$, and it does not admit any subsequence converging in $L^p(I)$.  Moreover, in case $p\neq 2$ there is a uniform constant $C=C(p)$ so that
\begin{eqnarray}
\int_0^1 \left|\sqrt{q(t)}~v'_j(t)\right|^p~dt&\leq& C(p)n^2\left[\int_{\frac{1}{n}+\frac{1}{n^{2+p}}}^\frac{2}{n}  t^{p^2} \left(t-\frac{1}{n}\right)^{1-p}~dt \right.\nonumber\\
&&\quad\quad\quad\;\;+\left. \int_{\frac{2}{n}}^{\frac{3}{n}-\frac{1}{n^{p+2}}}t^{p^2}\left(\frac{3}{n}-t\right)^{1-p}~dt\right]\nonumber\\
&\leq& C(p)n^{2-p^2}\left[\left(t-\frac{1}{n}\right)^{2-p}\bigg|_{\frac{1}{n}+\frac{1}{n^{2+p}}}^\frac{2}{n} - \left(\frac{3}{n}-t\right)^{2-p}\bigg|_\frac{2}{n}^{\frac{3}{n}-\frac{1}{n^{p+2}}}\right] \nonumber\\
&=&C(p)n^{2-p^2}\left[\frac{1}{n^{2-p}}-\frac{1}{n^{4-p^2}}\right]\nonumber\\
&\leq&C(p)<\infty\textrm{ for every $j$.}\nonumber
\end{eqnarray}
since $p>1$.  One can easily check the same for the case $p=2$ using logarithms.  Thus, the sequence is bounded in the one dimensional norm of $qH^{1,p}(I)$.  Therefore, $qH^{1,p}_0(I)\subset qH^{1,p}(I)$ are continuously but not compactly embedded in $L^p(I)$. 

In the $n$-dimensional case a similar result holds when one chooses $\Omega = I^n$,
$$Q(x) = \textrm{Diag}[x_1^{2p},1,...,1],$$
$v\equiv1$, $d\mu=dx$ and the sequence of functions 
$$u_j(x) = u_j(x_1,\hat{x}) = v_j(x_1)\psi(\hat{x})$$
for any $\psi\in C^\infty_0((0,1)^{n-1})$ with $\hat{x}=(x_2,...,x_n)$. It is important to note here that no global Sobolev inequality of order $p$ and gain $\sigma>1$ can hold on $I^n$ for the matrix $Q$.


\end{document}